\documentclass[12pt]{article}
\usepackage{amsfonts,amsmath,amssymb, amsthm,latexsym}
\usepackage[cp1251]{inputenc}
\usepackage[english, russian]{babel}
\usepackage[dvips]{graphicx}
\usepackage[dvips]{epsfig}
\usepackage{oldlfont}
\usepackage{srcltx}
\usepackage{amscd}
\usepackage{array,longtable}
\begin{document}

\begin{center}
\textbf{ O.Y. Kushel \\[0.5cm]
SPECTRAL PROPERTIES OF ONE CLASS \\ OF SIGN-SYMMETRIC MATRICES} \\[0.5cm]
\end{center}
\medskip

\begin{center}
{\bf Abstract.}
\end{center}

A $n \times n$ matrix $\mathbf{A}$, which has a certain
sign-symmetric structure ($\mathcal{J}$--sign-symmetric), is
studied in this paper. It's shown, that such a matrix is similar
to a nonnegative matrix. The existence of the second in modulus
positive eigenvalue $\lambda_2$ of a $\mathcal{J}$--sign-symmetric
matrix $\mathbf{A}$, or an odd number $k$ of simple eigenvalues,
which coincide with the $k$th roots of $\rho(A)^k$, is proved
under the additional condition, that its second compound matrix is
also $\mathcal{J}$--sign-symmetric. The conditions, when a
$\mathcal{J}$--sign-symmetric matrix with a
$\mathcal{J}$--sign-symmetric second compound matrix has complex
eigenvalues, which are equal in modulus to $\rho(A)$, are given.

\medskip

{\it Keywords: Totally positive matrices, Sign-symmetric matrices,
Nonnegative \linebreak matrices, Compound matrices, Exterior
powers, Gantmacher--Krein theorem, \linebreak eigenvalues.}

\medskip

{\it 2000 Mathematics Subject Classification: Primary 15A48,
Secondary 15A18, \linebreak 15A75.}

\medskip

\section{Introduction}

We first state a classical result of F.R. Gantmacher and M.G.
Krein (see, for example, [1]).

{\bf Theorem A (Gantmacher--Krein).} {\it If the matrix
$\mathbf{A}$ of a linear operator $A: {\mathbb{R}}^{n} \rightarrow
{\mathbb{R}}^{n}$ is positive together with its $j$th compound
matrices ${\mathbf A}^{(j)}$ $(1 < j \leq k)$ up to the order $k$,
then the operator $A$ has $k$ positive simple eigenvalues $0 <
\lambda_k < \ldots < \lambda_2 < \lambda_1$, with the positive
eigenvector $x_1$ corresponding to the maximal eigenvalue
$\lambda_1$, and the eigenvector $x_j$, which has exactly $j - 1$
changes of sign, corresponding to the $j$-th eigenvalue
$\lambda_j$.}

Let us remember, that a matrix, which satisfies the conditions of
theorem A, is called strictly $k$--totally positive ($STP_k$). If
a matrix $\mathbf{A}$ is nonnegative together with its $j$th
compound matrices ${\mathbf A}^{(j)}$ $(1 < j \leq k)$ up to the
order $k$, then $\mathbf{A}$ is called $k$--totally positive
($TP_k$). The following statement (see [1], p. 317, theorem 13)
easily comes out from the reasons of continuity: {\it if the
matrix $\mathbf{A}$ of a linear operator $A: {\mathbb{R}}^{n}
\rightarrow {\mathbb{R}}^{n}$ is $k$--totally positive, then the
operator $A$ has $k$ nonnegative eigenvalues $0 \leq \lambda_k
\leq \ldots \leq \lambda_2 \leq \lambda_1$.}

We next formulate the statement, which follows from the
infinite-dimensional statements, proved in [2]. This statement
generalizes the concept of total positivity.

{\bf Theorem B.} {\it Let the matrix ${\mathbf A}$ of a linear
operator $A: {\mathbb R}^n \rightarrow {\mathbb R}^n$ be similar
to a nonnegative irreducible matrix. Let its second compound
matrix ${\mathbf A}^{(2)}$ be also similar to a nonnegative
irreducible matrix. Then one of the following two cases takes
place:
\begin{enumerate}
\item[\rm (1)] The operator $A$ has two positive
simple eigenvalues $\lambda_{1}$, $\lambda_{2}$:
$$\rho(A) = \lambda_{1} > \lambda_{2} > 0. $$
Moreover, if the operator $A$ has $l > 1$ eigenvalues
$\lambda_{2}, \lambda_{3}, \ldots, \lambda_{l + 1} $, equal in
modulus to $\lambda_{2}$, then each of them is simple, and they
coincide with the $l$th roots from $\lambda_{2}^l$.
\item[\rm (2)] There is just three eigenvalues on the
spectral circle $|\lambda| = \rho(A)$. Each of them is simple, and
they coincide with the $3$th roots of $(\rho(A))^3$.
\end{enumerate}}

The proof of theorem B is based on the well-known Perron-Frobenius
statement, that a nonnegative irreducible matrix has the maximal
in modulus simple positive eigenvalue, and very simple reasoning,
that similarity transformations do not change the spectrum of a
matrix. Let us note, that the peripheral spectrum of the operator
$A$ not always will be real. This is the principal difference
between the statement of theorem B and the statements, proved in
[1]. The necessary and sufficient condition of the existence of
complex eigenvalues on the spectral circle $|\lambda| = \rho(A)$
in the case, when both the matrices ${\mathbf A}$ and ${\mathbf
A}^{(2)}$ are similar to nonnegative irreducible matrices with
diagonal matrices of similarity transformations, was obtained in
[3].

The class of matrices, which are similar to nonnegative matrices
together with their second compound matrices is studied in this
paper (later on we shall not assume the additional condition of
irreducibility). The problem of the detailed description of the
spectrum of such matrices is examined. The obtained results are
based on the Frobenius theorems of the structure of reducible
matrices and combinatorial reasons. The main question, raised in
this paper, is: when will such matrices have complex eigenvalues
on the largest spectral circle?

\section{Tensor and exterior square of the space ${\mathbb{R}}^{n}$}

Let us briefly remind basic definitions and statements, related
with the tensor and exterior squares of the space
${\mathbb{R}}^{n}$ (for more detailed information see [4], [5],
[3]).

Later on, as the author thinks, it will be more convenient to
consider the space ${\mathbb R}^n$ as the space of real-valued
functions $x: \{1, \ \ldots, \ n\} \rightarrow {\mathbb R}$,
defined on the finite set of indices $\{1, \ \ldots, \ n \}$.
Denote such a space by ${\mathbb X}$. The basis in the space
${\mathbb X}$ consists of the functions $e_i$, for which $e_i(j) =
\delta_{ij}$.

The tensor square ${\mathbb X} \otimes {\mathbb X}$ of the space
${\mathbb X}$ is the space of all functions, defined on the set
$\{1, \ \ldots, \ n\} \times \{1, \ \ldots, \ n\}$. The space
${\mathbb X} \otimes {\mathbb X}$ is isomorphic to the space
${\mathbb R}^{n^2}$.

 The exterior square ${\mathbb X}\wedge
{\mathbb X}$ of the space ${\mathbb X}$ is the space of all
antisymmetric functions (i.e. functions $f(i,j)$, for which the
equality $f(i,j) = - f(j,i)$ is true), defined on the set $\{1, \
\ldots, \ n\} \times \{1, \ \ldots, \ n\}$. It is known, that
${\mathbb X}\wedge {\mathbb X}$ coincides with the linear span of
all exterior products $x \wedge y \ \ (x,y \in {\mathbb X})$,
which acts according to the rule:
$$(x \wedge y)(i,j) =
x(i)y(j) - x(j)y(i).$$ The space ${\mathbb X}\wedge {\mathbb X}$
is isomorphic to the space ${\mathbb X} (W\setminus\Delta)$ of
real-valued functions, defined on the set $W\setminus\Delta$,
where $W$ is a subset of $\{1, \ \ldots, \ n\} \times \{1, \
\ldots, \ n\}$, which satisfies the following conditions:$$ W\cup
\widetilde{W} = (\{1, \ \ldots, \ n\} \times \{1, \ \ldots, \
n\}); \eqno(1)$$
$$W \cap \widetilde{W} = \Delta. \eqno(2)$$ (Here $\widetilde{W}
= \{(j,i): \ (i,j) \in W\}$; $\Delta = \{(i,i): \ i =1, \ \ldots,
\ n\}$).

The following equality is true for the power of the set
$W\setminus\Delta$:
$$N(W\setminus\Delta) = \frac{n^2 - n}{2} = C_n^2.$$

It comes out from the above reasoning, that the space ${\mathbb X}
(W\setminus\Delta)$ is isomorphic to the space ${\mathbb
R}^{C_n^2}$.

Every subset $W \subset \{1, \ \ldots, \ n\} \times \{1, \ \ldots,
\ n\}$ defines a binary relation on the set $\{1, \ \ldots, \ n\}$
(see, for example, [6]). If the set $W$ satisfies conditions (1)
and (2), and, in addition, the inclusion $(i,k) \in W$ follows
from the inclusions $(i,j) \in W$ and $(j,k) \in W$ for any $i, j,
k \in \{1, \ \ldots, \ n\}$ (i.e. the set $W$ possesses the
property of {\it transitivity}), then the relation, defined by the
set $W$, is a linear order relation.

 Every subset $W \subset \{1, \ \ldots, \ n\} \times \{1, \ \ldots,
\ n\}$, which satisfies conditions (1) and (2), uniquely defines a
basis $\{e_i \wedge e_j\}_{(i,j)\in W \setminus \Delta}$ in
${\mathbb R}^n \wedge {\mathbb R}^n$, which consists of the
exterior products of the initial basic vectors. Such a basis in
${\mathbb X}\wedge {\mathbb X}$, constructed with respect to the
set $W$, is called {\it a $W$--basis}.

\section{The second compound matrix and the exterior square
of a linear operator}

The exterior square $A \wedge A$ of the operator $A: {\mathbb
X}\rightarrow {\mathbb X}$ acts in the space ${\mathbb
X}\wedge{\mathbb X}$ according to the rule: $$ (A \wedge A)(x
\wedge y) = Ax \wedge Ay. $$ Let the operator $A$ be defined by
the matrix ${\mathbf A} = \{a_{ij}\}_{i,j = 1}^n$ in the basis
$\{e_i\}_{i = 1}^n$. Then the matrix of the operator $A \wedge A$
in the $W$--basis $\{e_i \wedge e_j\}_{(i,j)\in W \setminus
\Delta}$ coincides with the $W$--matrix $ {\mathbf A}_W^{(2)}$.
(Here the matrix $ {\mathbf A}_W^{(2)}$ consists of minors
$A\begin{pmatrix}
  i & j \\
  k & l
\end{pmatrix}$ of the initial matrix $\mathbf{A}$, formed of the rows with
numbers $i$ and $j$ and the columns with numbers $k$ and $l$,
where $ (i,j), \ (k,l) \in (W \setminus \Delta)$. The minors
$A\begin{pmatrix}
  i & j \\
  k & l
\end{pmatrix}$ are
numerated in the lexicographic order.)

Note, that if the set $W$ coincides with the set $M = \{(i,j) \in
\{1, \ \ldots, \ n\} \times \{1, \ \ldots, \ n\}: \ i \leq j \}$,
then the corresponding $W$--basis is $\{e_i \wedge e_j\}_{1 \leq
i<j \leq n}$, i.e. the canonical basis in the space ${\mathbb R}^n
\wedge {\mathbb R}^n$, and the corresponding $W$--matrix is a
matrix, which consists of minors $A\begin{pmatrix}
  i & j \\
  k & l
\end{pmatrix}$, where $1\leq i<j \leq n, \ 1\leq k < l \leq n$,
i.e. the second compound matrix of the matrix $\mathbf{A}$.

The following theorem about the eigenvalues of a $W$--matrix is
true (see [3]).

{\bf Theorem 1.} \textit{Let $W$ be a subset of the set $\{1, \
\ldots, \ n\} \times \{1, \ \ldots, \ n\}$, which satisfies
conditions (1) and (2). Let $\{\lambda_{i}\}_{i = 1}^n$ be the set
of all eigenvalues of the matrix ${\mathbf A}$, repeated according
to multiplicity. Then all the possible products of the type
$\{\lambda_{i} \lambda_{j} \}$, where $1 \leq i < j \leq n$, forms
the set of all the possible eigenvalues of the corresponding
$W$--matrix $ {\mathbf A}_W^{(2)}$, repeated according to
multiplicity}.

In the case $W = M$ theorem 1 turns into the Kronecker theorem
(see [1], p. 80, theorem 23) about the eigenvalues of the second
compound matrix.

\section{Operators, which leave invariant a cone in ${\mathbb R}^n$, and their matrices}

Let us remember some widely used definitions. An operator $A:
{\mathbb R}^n \rightarrow {\mathbb R}^n$ is called {\it
nonnegative}, if it leaves invariant the cone of nonnegative
vectors in ${\mathbb R}^n$. It's well known, that the operator $A$
is nonnegative if and only if its matrix $\mathbf{A}$ is
nonnegative, i.e. all the elements $a_{ij} \ (i,j = 1, \ \ldots, \
n)$ of the matrix $\mathbf{A}$ are nonnegative.

The following theorem is perhaps the best known part of the theory
of nonnegative operators (see, for example, [7], p. 14, theorem
4.2).

{\bf Theorem 2 (Perron).} {\it Let the matrix $\mathbf A$ of a
linear operator $A:{\mathbb{R}}^{n} \rightarrow {\mathbb{R}}^{n}$
be nonnegative. Then the spectral radius $\rho(A) \geq 0$ is a
nonnegative eigenvalue of the operator $A$, with the nonnegative
eigenvector $x_1$ corresponding to it.}

However, this result is also correct for any matrix, similar to a
nonnegative matrix (since a similarity transformation preserves
the spectrum of a matrix). The question, how can we see if an
arbitrary matrix is similar to a nonnegative matrix, was raised in
[3]. The answer was given for the special case, when the matrix of
the similarity transformation is diagonal. The following
definition was first introduced in [3].

A matrix ${\mathbf A}$ of a linear operator $A: {\mathbb R}^n
\rightarrow {\mathbb R}^n$ is called {\it ${\mathcal
J}$--sign-symmetric}, if there exists such a subset ${\mathcal J}
\subseteq \{1, \ \ldots, \ n\}$, that both the conditions (a) and
(b) are true:
\begin{enumerate}
\item[\rm (a)] the inequality $a_{ij} \leq 0$ follows from the inclusions $i
\in {\mathcal J}$, $j \in \{1, \ \ldots, \ n\}\setminus {\mathcal
J}$ and from the inclusions $j \in {\mathcal J}$, $i \in \{1, \
\ldots, \ n\}\setminus {\mathcal J}$ for any two numbers $i,j$;

\item[\rm (b)] one of the inclusions $i \in {\mathcal J}$, $j \in \{1, \
\ldots, \ n\}\setminus {\mathcal J}$ or $j \in {\mathcal J}$, $i
\in \{1, \ \ldots, \ n\}\setminus {\mathcal J}$ follows from the
strict inequality $a_{ij} < 0$.
\end{enumerate}

Note, that the subset ${\mathcal J}$ in the definition of
${\mathcal J}$--sign-symmetricity in the case of an arbitrary
matrix ${\mathbf A}$ is not uniquely defined.

The following statement of the similarity of a ${\mathcal
J}$--sign-symmetric matrix to a nonnegative matrix was proved in
[3] (see [3], theorem 4).

{\bf Theorem 3.} {\it Let ${\mathbf A}$ be a ${\mathcal
J}$--sign-symmetric matrix. Then it can be represented in the
following form:
$${\mathbf A} = {\mathbf D} \widetilde{{\mathbf A}} {\mathbf
D}^{-1}, \eqno(3)$$ where $\widetilde{{\mathbf A}}$ is a
nonnegative matrix, ${\mathbf D}$ is a diagonal matrix, which
diagonal elements are equal to $\pm 1$.}

{\bf Corollary.} {\it Let the matrix $\mathbf A$ of a linear
operator $A:{\mathbb{R}}^{n} \rightarrow {\mathbb{R}}^{n}$ be
${\mathcal J}$--sign-symmetric. Then the spectral radius $\rho(A)
\geq 0$ is a nonnegative eigenvalue of the operator $A$.}

This statement has the following geometric meaning: if the matrix
$\mathbf A$ of a linear operator $A$ is ${\mathcal
J}$--sign-symmetric in some basis $\{e_i\}_{i = 1}^n$, then the
operator $A$ leaves invariant one of the cones, spanned on the
vectors $\{ \pm e_i\}_{i = 1}^n$.

Later on operators with positive spectral radius will play an
important role. Let us formulate a sufficient criteria of the
positivity of $\rho(A)$.

{\bf Lemma 1.} {\it Let the matrix $\mathbf A$ of a linear
operator $A:{\mathbb{R}}^{n} \rightarrow {\mathbb{R}}^{n}$ be
${\mathcal J}$--sign-symmetric. Let at least one element $a_{ii}$,
situated on the principal diagonal of $\mathbf A$, be not equal to
zero. Then the spectral radius $\rho(A)
> 0$ is a positive eigenvalue of the operator $A$.}

$\square$ It follows from the definition of ${\mathcal
J}$--sign-symmetricity, that all the elements $a_{ii} \ (i= 1, \
\ldots, \ n)$ are nonnegative. Since at least one of $a_{ii}$ is
non-zero, then the sum $ \sum_{i = 1}^n a_{ii}$ is positive. It's
clear (see, for example, [5]), that:
 $$tr(A) = \sum_{i = 1}^n \lambda_i = \sum_{i = 1}^n a_{ii}.$$
(Here $\{\lambda_i\}_{i = 1}^n $ is the set of all eigenvalues of
the operator $A$, repeated according to multiplicity.) So we have,
that $$\rho(A) \geq \frac{1}{n}\sum_{i = 1}^n \lambda_i =
\frac{1}{n}\sum_{i = 1}^n a_{ii}
> 0.$$ $\blacksquare$

This criteria is sufficient, but, as it'll be shown later, not
necessary for the positivity of $\rho(A)$.

Now we'll require some notations. Given the set $\alpha$ of $k \ \
(1 \leq k \leq n)$ indices $i_1, \ \ldots, \ i_k  \ \ (1 \leq i_1
< \ldots
< i_k \leq n)$. Then $${\mathbf A}(\alpha) = {\mathbf A}\left[%
\begin{array}{ccc}
  i_1 & \ldots & i_k \\
  i_1 & \ldots & i_k \\
\end{array}%
\right] $$ denotes a $k \times k$ submatrix of the matrix
$\mathbf{A}$, formed of the rows with numbers $i_1, \ \ldots, \
i_k$ and the columns with the same numbers. The submatrix
${\mathbf A}(\alpha)$ is called a {\it principal submatrix} of the
matrix $\mathbf{A}$ and it represents the matrix of the
restriction of the operator $A$ to a $k$-dimensional coordinate
subspace, spanned on the basic vectors $e_{i_1}, \ \ldots, \
e_{i_k}$.

Any submatrix of a nonnegative matrix is obviously nonnegative.
Let us prove the corresponding property of principal submatrices
of a ${\mathcal J}$--sign-symmetric matrix.

{\bf Lemma 2.} {\it Let the matrix ${\mathbf A}$ of a linear
operator $A$ be ${\mathcal J}$--sign-symmetric. Then every
principal submatrix of the matrix ${\mathbf A}$ is also ${\mathcal
J}$--sign-symmetric.}

$\square$ Let ${\mathbf A}(\alpha)$ is an arbitrary principal
submatrix of the matrix ${\mathbf A}$. It's obvious, that the set
${\mathcal J}\cap \alpha$ satisfy both the conditions (a) and (b)
in the definition of ${\mathcal J}$--sign-symmetricity of the
matrix ${\mathbf A}(\alpha)$. $\blacksquare$

\section{Irreducible operators and their matrices}

A linear operator $ A: {\mathbb R}^n \rightarrow {\mathbb R}^n$ is
called {\it irreducible}, if it has no invariant $k$-dimensional
coordinate subspaces with $0 < k < n$. Otherwise it's called {\it
reducible}. An operator $A$ is irreducible (reducible) if and only
if its matrix $\mathbf{A}$ is irreducible (respectively
reducible). Remember, that a matrix ${\mathbf A}$ is called {\it
reducible}, if there exists a permutation of coordinates such
that: $${\mathbf P}^{-1}{\mathbf A}{\mathbf P} =
\begin{pmatrix}
  {\mathbf A}_1 & 0 \\
  {\mathbf B} & {\mathbf A}_2
\end{pmatrix}, \eqno(4)$$ where ${\mathbf P}$ is an $n \times n$ permutation matrix
(each row and each column have exactly one 1 entry and all others
0), ${\mathbf A}_1$, ${\mathbf A}_2$ are square matrices.
Otherwise the matrix $\mathbf A$ is called irreducible.

We next formulate some important properties of irreducible
operators (see, for example, [7], [8]).

{\it The following properties of a linear operator $A: {\mathbb
R}^n \rightarrow {\mathbb R}^n$ are equivalent:

(i) Its matrix $\mathbf{A}$ is irreducible;

(ii) Its matrix $\mathbf{A}$ has a "path of irreducibility", i.e.
such a set of indices $\{j_0, j_1, \ldots , j_s\}$, that $j_0
\not= j_1, \ j_1 \not= j_2, \ldots, j_{s-1} \not= j_s$, every
index of $\{1, \ldots , n\}$ coincides with one of the indices
$\{j_0, j_1, \ldots , j_s \}$ and $$a_{j_0 \, j_1}, \ a_{j_1 \,
j_2}, \ \ldots \ a_{j_{s-1} \, j_s}, \ a_{j_s \, j_0} > 0.$$

(iii) Every nonnegative eigenvector of the operator $A$ is
positive.}

An irreducible matrix $\mathbf{A}$ is called {\it imprimitive}, if
there exists a permutation of coordinates such that: $${\mathbf
P}^{-1}{\mathbf A}{\mathbf P} =
\begin{pmatrix}
  0 & {\mathbf A}_{12} & 0 & \ldots & 0 \\
  0 & 0 & {\mathbf A}_{23} & \ldots & 0 \\
  \ldots & \ldots & \ldots & \ldots & \ldots \\
  0 & 0 & 0 & \ldots & {\mathbf A}_{h-1 h} \\
  {\mathbf A}_{h1} & 0 & 0 & \ldots & 0
\end{pmatrix}, \eqno(5)$$ where ${\mathbf P}$ is an $n \times n$ permutation matrix
, ${\mathbf A}_{i \ i+1}$ $(i = 1, \ \ldots, \ h-1)$ and ${\mathbf
A}_{h1}$ are square matrices.
 Otherwise the matrix $\mathbf{A}$ is called {\it primitive}.
An operator $A$ is imprimitive (primitive) if and only if its
matrix is imprimitive (primitive).

The following sufficient criteria of primitivity was proved in [7]
(see [7], p. 49, corollary 1.1): {\it if the matrix $\mathbf A$ of
the operator $A$ is irreducible, and $tr(A)
> 0$, then $\mathbf A$ is primitive.}

Let us remember the widely known Frobenius theorem on the spectrum
of irreducible operators.

{\bf Theorem 4 (Frobenius).} {\it Let the matrix $\mathbf A$ of a
linear operator $A$ be nonnegative and irreducible. Then the
spectral radius $\rho(A)
> 0$ is a simple positive eigenvalue of the operator $A$,
with the corresponding positive eigenvector $x_1$. If $h$ is a
number of the eigenvalues of the operator $A$, which are equal in
modulus to $\rho(A)$, then all of them are simple and they
coincide with the $h$-th roots of $(\rho(A))^h$. Moreover, the
spectrum of the operator $A$ is invariant under rotations by
$\frac{2\pi}{h}$ about the origin.}

The number of the eigenvalues, which are equal in modulus to
$\rho(A)$ is called {\it the index of imprimitivity} of an
irreducible operator $A$ and denoted $h(A)$. If $h(A) > 1$, then
it is equal to the number $h$ of the square blocs in the form (5)
of the matrix $\mathbf{A}$. If $\mathbf A$ is primitive, then
$h(A) = 1$.

It's easy to see, that the similarity transformation (3) with a
diagonal matrix ${\mathbf D}$ preserves the property of
irreducibility and the index of imprimitivity of the operator $A$
as well, as other spectral properties (see [3]). So we have:

{\bf Theorem 5.} {\it Let the matrix $\mathbf A$ of a linear
operator $A$ be ${\mathcal J}$--sign-symmetric and irreducible.
Then the spectral radius $\rho(A)
> 0$ is a simple positive eigenvalue of the operator $A$. If $h$ is a
number of the eigenvalues of the operator $A$, which are equal in
modulus to $\rho(A)$, then all of them are simple and they
coincide with the $h$-th roots of $(\rho(A))^h$. Moreover, the
spectrum of the operator $A$ is invariant under rotations by
$\frac{2\pi}{h}$ about the origin.}

And the following criteria is true.

{\bf Lemma 3.} {\it Let the matrix $\mathbf A$ of the operator $A$
be ${\mathcal J}$--sign-symmetric and \linebreak irreducible. If
at least one element $a_{ii}$, situated on the principal diagonal
of $\mathbf A$, is not equal to zero, then $\mathbf A$ is
primitive. In turn, if $\mathbf A$ is imrimitive, then its
principal diagonal contains only zeroes.}

 $\square$ The proof of lemma 3 comes out from the given above
sufficient criteria of the primitivity of nonnegative matrices.
$\blacksquare$

Note, that {\it if the matrix ${\mathbf A}$ is ${\mathcal
J}$--sign-symmetric and irreducible, then the set ${\mathcal J}$
in the definition of ${\mathcal J}$--sign-symmetricity is uniquely
defined (up to the set $\{1, \ \ldots, \ n\}\setminus {\mathcal
J}$)}.

\section{Reducible operators and their matrices}

The statement, which helps to bring the study of reducible
operators in the finite-dimensional space ${\mathbb R}^n$ to the
study of irreducible operators is well-known (see, for example,
[8]).

{\bf Theorem 6 (Frobenius).} {\it Let $A$ be a nonnegative
reducible matrix. Then there exists a permutation of coordinates
such that:
$${\mathbf P}^{-1}{\mathbf A}{\mathbf P} = \widehat{{\mathbf A}},
$$ where ${\mathbf P}$ is an $n \times n$ permutation matrix,
$\widehat{{\mathbf A}}$ is a block-triangular form, where the
finite number $l \leq n$ of square irreducible blocs $A_j \ \ (j =
1, \ldots , l)$ are situated on the principal
 diagonal, and zero elements are situated above the principal
diagonal:
 $$ A = \left(\begin{array}{cccccccc}A_1 &  0  & \ldots &  0  & 0 & 0 & \ldots & 0 \cr
 0  & A_2 & \ldots &  0  & 0 & 0 & \ldots & 0 \cr
 \ldots & \ldots & \ldots & \ldots & \ldots & \ldots & \ldots & \ldots \cr
 0  &  0  & \ldots & A_r & 0 & 0 & \ldots & 0 \cr
 B_{r+1 \, 1} & B_{r+1 \, 2} & \ldots & B_{r+1 \, r} & A_{r+1} & 0 & \ldots & 0 \cr
 B_{r+2 \, 1} & B_{r+2 \, 2} & \ldots & B_{r+2 \, r} &
 B_{r+2 \, r+1} & A_{r+2} & \ldots & 0 \cr
 \ldots & \ldots & \ldots & \ldots & \ldots & \ldots & \ldots & \ldots \cr
 B_{l \, 1} & B_{l \, 2} & \ldots & B_{l \, r} & B_{l \, r+1} &
 B_{l \, r+2} & \ldots & A_l \cr\end{array}\right) \eqno(6)$$
Such a representation is uniquely defined (up to the numeration of
the blocs).

 The spectral radius $\rho(A)$ is an eigenvalue of the
operator $A$ with the nonnegative eigenvector $x_1$ corresponding
to it. Moreover, the following equalities are true:
$$ \sigma_p(A) = \bigcup_{j=1}^l \sigma_p(A_j), \quad \rho(A) =
 \max_{j = 1, \ldots , l} \  \{\rho(A_j)\},$$
where $\sigma_p(A_j)$ are the point spectra (i.e. the sets of all
eigenvalues), and $\rho(A_j)$ are the spectral radii of the
irreducible blocs $A_j \ \ (j = 1, \ldots , l)$.}

It's easy to see, that {\it if the matrix ${\mathbf A}$ is
${\mathcal J}$--sign-symmetric and reducible, then there is $2^l$
 possible ways of constructing the set ${\mathcal J}$ in the definition of ${\mathcal
J}$--sign-symmetricity}.

\section{The connection between the sets ${\mathcal J}$ and $W$.}

Later on we'll impose the condition of ${\mathcal
J}$--sign-symmetricity to the second compound matrices ${\mathbf
A}^{(2)}$. So we are interested in some special properties of
${\mathbf A}^{(2)}$. The following statement about the link
between
 the structure of the second compound matrices ${\mathbf
A}^{(2)}$ and the structure of the matrix ${\mathbf A}_W^{(2)}$,
constructed with respect to a set $W \in \{1, \ \ldots, \
n\}\times \{1, \ \ldots, \ n\}$, which satisfies conditions (1)
and (2), is true.

 {\it Let the second compound matrix ${\mathbf
A}^{(2)}$ of a matrix ${\mathbf A}$ be ${\mathcal
J}$--sign-symmetric. Then there exists such a set $W \in \{1, \
\ldots, \ n\}\times \{1, \ \ldots, \ n\}$, which satisfies
conditions (1) and (2), that the corresponding $W$--matrix
${\mathbf A}_W^{(2)}$ is nonnegative. Moreover, if ${\mathbf
A}^{(2)}$ is irreducible, then ${\mathbf A}_W^{(2)}$ is also
irreducible} (see [3], theorem 6).

This statement easily follows from the theorem 3 and the given
above fact, that the matrix of the exterior square $A \wedge A$ of
the operator $A$ in the $W$--basis \linebreak $\{e_i \wedge
e_j\}_{(i,j)\in W \setminus \Delta}$ coincides with the
$W$--matrix $ {\mathbf A}_W^{(2)}$.

The method of constructing the set $W$, for which the
corresponding $W$--matrix ${\mathbf A}_W^{(2)}$ is positive, by
the set ${\mathcal J}$ in the definition of ${\mathcal
J}$-sign-symmetricity of $ {\mathbf A}^{(2)}$ is given in [3] (see
the proof of theorems 5 and 6):

The pair $(i,j) \in W$ if and only if one of the following two
cases takes place:
\begin{enumerate}
\item[\rm (a)] $i < j$, and the number $\alpha$ of the pair $(i,j)$
(in the lexicographic numeration), belongs to the set ${\mathcal
J}$;

\item[\rm (b)] $i > j$, and the number $\widetilde {\alpha}$
of the pair $(j,i)$ belongs to the set $\{1, \ \ldots, \
C_n^2\}\setminus {\mathcal J}$.\end{enumerate}

It's easy to see, that the number of different ways of the
constructing of the set $W$ is equal to the number of different
ways of the constructing of the set ${\mathcal J}$ in the
definition of ${\mathcal J}$-sign-symmetricity of the second
compound matrix.

 Let us generalize the given above method of constructing the
set $W$ to the case of ${\mathcal J}$--sign-symmetric matrix. Let
${\mathbf A}$ be a ${\mathcal J}$--sign-symmetric matrix, and let
${\mathcal J}$ be a subset of the set $\{1, \ \ldots, \ n\}$ in
the definition of ${\mathcal J}$--sign-symmetricity (i.e. such a
subset, that the conditions (a) and (b) are true). Let ${\mathbf
A}^{(2)}$ be a ${\mathcal J}$--sign-symmetric matrix. Let
$\widetilde{{\mathcal J}}$ be a subset of $\{1, \ \ldots, C_n^2\}$
in the definition of ${\mathcal J}$--sign-symmetricity for the
matrix ${\mathbf A}^{(2)}$. Let us construct a set
$\widehat{W}({\mathcal J}, \widetilde{{\mathcal J}}) \subseteq
(\{1, \ \ldots, \ n\} \times \{1, \ \ldots, \ n\})$ with respect
to the sets ${\mathcal J}$ and $\widetilde{{\mathcal J}}$ by the
following way.

A pair $(i,j)$ belongs to the set $\widehat{W}({\mathcal J},
\widetilde{{\mathcal J}})$ if and only if one of the following
four cases takes place:
\begin{enumerate}
\item[\rm (a)] $i < j$, both the numbers $i,j$ belong either to the set ${\mathcal
J}$, or to the set $\{1, \ \ldots, \ n\}\setminus {\mathcal J}$,
and the number $\alpha$, corresponding to the pair $(i,j)$ in the
lexicographic numeration, belongs to the set $\widetilde{{\mathcal
J}}$;

\item[\rm (b)] $i < j$, one of the numbers $i,j$ belongs to the set ${\mathcal
J}$, and the other belongs to the set $\{1, \ \ldots, \
n\}\setminus {\mathcal J}$, and the number $\alpha$, corresponding
to the pair $(i,j)$ in the lexicographic numeration, belongs to
the set $\{1, \ \ldots, \ C_n^2\}\setminus \widetilde{{\mathcal
J}}$;

\item[\rm (c)] $i > j$, both the numbers $i,j$ belong either to the set ${\mathcal
J}$, or to the set $\{1, \ \ldots, \ n\}\setminus {\mathcal J}$,
and the number $\alpha$, corresponding to the pair $(j, i)$ in the
lexicographic numeration, belongs to the set $\{1, \ \ldots, \
C_n^2\}\setminus \widetilde{{\mathcal J}}$;

\item[\rm (d)]  $i > j$, one of the numbers $i,j$ belongs to the set ${\mathcal
J}$, the other belongs to the set $\{1, \ \ldots, \ n\}\setminus
{\mathcal J}$, and the number $\alpha$, corresponding to the pair
$(j, i)$ in the lexicographic numeration, belongs to the set
$\widetilde{{\mathcal J}}$.

\end{enumerate}

 As it was noticed above, the set $W({\mathcal J}, \widetilde{{\mathcal J}})$
is called {\it transitive}, if the inclusion $(i,k) \in
W({\mathcal J}, \widetilde{{\mathcal J}})$ follows from the
inclusions $(i,j) \in W({\mathcal J}, \widetilde{{\mathcal J}})$
and $(j,k) \in W({\mathcal J}, \widetilde{{\mathcal J}})$ for any
$i, j, k \in \{1, \ \ldots, \ n\}$

The set $W({\mathcal J}, \widetilde{{\mathcal J}})$ is obviously
not uniquely defined, but there is a finite number of the possible
ways of its constructing.

\section{Generalization of the Gantmacher--Krein theorems to the
case of an irreducible $2$--totally ${\mathcal J}$--sign-symmetric
matrix.}

Now we'll prove the theorem of the spectral properties of an
irreducible ${\mathcal J}$--sign-symmetric matrix with a
${\mathcal J}$--sign-symmetric second compound matrix, using the
following statement, proved in [3] (see [3], theorem 12).

 {\bf Theorem 7.} {\it Let ${\mathbf A}$ be a ${\mathcal
J}$--sign-symmetric matrix. Let its second compound matrix
${\mathbf A}^{(2)}$ be also ${\mathcal J}$--sign-symmetric. Let
one of the sets $\widehat{W}({\mathcal J}, \widetilde{{\mathcal
J}})$ be transitive. Then the two largest in modulus eigenvalues
of the operator $A$ are nonnegative.}

The following statement about the spectrum of a ${\mathcal
J}$--sign-symmetric irreducible matrix with a ${\mathcal
J}$--sign-symmetric second compound matrix is true.

{\bf Theorem 8.} {\it Let the matrix ${\mathbf A}$ of a linear
operator $A$ be ${\mathcal J}$--sign-symmetric and irreducible.
Let its second compound matrix ${\mathbf A}^{(2)}$ be ${\mathcal
J}$--sign-symmetric. Then the operator $A$ has a simple positive
eigenvalue $\lambda_{1} = \rho(A)$. Moreover, one of the following
two cases takes place:
\begin{enumerate}
 \item[\rm (1)] If at least one of the possible sets $W({\mathcal J},
\widetilde{{\mathcal J}})$ is transitive, then $h(A) = 1$, the
second in modulus eigenvalue $\lambda_{2}$ of the operator $A$ is
nonnegative and different in modulus from the first eigenvalue:
 $$0 \leq \lambda_2 < \lambda_{1}.$$

 \item[\rm (2)] If all the possible sets $W({\mathcal J},
\widetilde{{\mathcal J}})$ are not transitive, then there is an
odd number $k \geq 1$ of eigenvalues on the spectral circle
$|\lambda| = \rho(A)$. All of them are simple and coincide with
the $k$th roots of $(\rho(A))^k$.
\end{enumerate}}

$\square$ Suppose a transitive set $W({\mathcal J},
\widetilde{{\mathcal J}})$ exists. Enumerate the eigenvalues of
the matrix ${\mathbf A}$ (repeated according to multiplicity) in
order of decrease of their modules:
$$|\lambda_{1}| \geq | \lambda_{2}| \geq |\lambda_{3}| \geq \ldots
\geq |\lambda_{n}|.$$ Then, applying theorem 7, we get, that the
two largest in modulus eigenvalues of the operator $A$ are
nonnegative: $0 \leq \lambda_2 \leq \lambda_1 = \rho(A)$. It
means, that either $\lambda_1$ is a multiple eigenvalue (this
contradicts the irreducibility of the operator $A$) or there is
only one eigenvalue on the spectral circle $|\lambda| = \rho(A)$.
I.e. the equality $h(A) = 1$ is true and the inequalities $0 \leq
\lambda_2 < \lambda_{1}$ are also true.

Let all the possible sets $W({\mathcal J}, \widetilde{{\mathcal
J}})$ be not transitive. Let $h(A) = k$ be the index of
imprimitivity of the operator $A$. Let us apply the Frobenius
theorem to the irreducible operator $A$. We get, that all the
eigenvalues of the operator $A$, equal in modulus to $\rho(A)$,
are simple and they coincide with the $k$th roots of $\rho(A)^k$.

Show, that the number $k$ is odd in this case. Suppose the
opposite: let $k$ be even. Examine the exterior square $A \wedge
A$ of the operator $A$ and its matrix ${\mathbf A}^{(2)}$. It
follows from theorem 1, that all the eigenvalues $A \wedge A$,
equal in modulus to $\rho(A \wedge A)$, form all the possible
products of the type $\lambda_j\lambda_m$, where $1\leq j<m \leq
k$, and $\lambda_1, \ \ldots, \ \lambda_k$ are all the eigenvalues
of the operator $A$, equal in modulus to $\rho(A)$. If $k = 2$,
then there is only one negative eigenvalue equal to $- \rho(A)^2$
on the spectral circle $|\lambda| = \rho(A \wedge A)$. This fact
contradicts the Perron theorem (applying the Perron theorem to the
operator $A \wedge A$, we get, that the spectral radius $\rho(A
\wedge A)$ is a nonnegative eigenvalue of $A \wedge A$). Now let
$k = 4, \ 6, \ 8, \ \ldots$. Examine $\lambda_{j} =
\rho(A)e^{\frac{2\pi(j-1)i}{k}} \ (j = 1, \ldots, k)$ --- all the
eigenvalues of the operator $A$, equal in modulus to $\rho(A)$. It
follows from theorem 1, that all the eigenvalues, equal in modulus
to $\rho(A \wedge A)$, can be written in the form
$\lambda_{j}\lambda_{m} =
\rho(A)^2e^{\frac{2\pi(j-1)i}{k}}e^{\frac{2\pi(m-1)i}{k}}, $ where
$1 \leq j<m \leq k$ (the general number of such eigenvalues is
$C_k^2 = \frac{k(k-1)}{2}$). It's easy to see, that there are
multiple eigenvalues among them. That is why the operator $A
\wedge A$ is reducible. It follows from the Frobenius theorem of
nonnegative matrices, that the number $s$ of irreducible blocs
${\mathcal A}_j$ with the property $\rho({\mathcal A}_j) = \rho(A
\wedge A)$ on the principal diagonal in representation (6) of the
reducible operator $A \wedge A$ is equal to the multiplicity of
the nonnegative eigenvalue $\rho(A \wedge A) = \rho(A)^{2}$. Let
us calculate this multiplicity. For this we write down all the
products of the form $\lambda_{j}\lambda_{m}$, where $1\leq j<m
\leq k$, equal to $\rho(A)^{2}$:
$$\rho(A)^{2} = \rho(A)^{2}e^{\frac{2\pi i}{k}}e^{\frac{2\pi(k-1)i}{k}} $$
$$\rho(A)^{2} = \rho(A)^{2}e^{\frac{4\pi i}{k}}e^{\frac{2\pi(k-2)i}{k}} $$
$$\rho(A)^{2} = \rho(A)^{2}e^{\frac{6\pi i}{k}}e^{\frac{2\pi(k-3)i}{k}} $$
$$\ldots \ \ \ \ldots \ \ \ \ldots $$
$$\rho(A)^{2} = \rho(A)^{2}e^{\frac{2\pi(\frac{k}{2}-1)i}{k}}
e^{\frac{2\pi(\frac{k}{2}+1)i}{k}} $$

It's easy to see, that the number of such products is equal to
$\frac{k}{2} - 1$. As it follows, the number $s$ of the
irreducible blocs ${\mathcal A}_j$ with $\rho({\mathcal A}_j) =
\rho(A \wedge A)$ is equal to $\frac{k}{2} - 1$. Since all the
eigenvalues, equal in modulus to $\rho(A\wedge A)$, are products
of the different $k$th roots of $\rho(A)^k$, they are $k$th roots
of $\rho(A)^{2k}$. As it follows, the index of imprimitivity
$h({\mathcal A}_j)$ can not be greater than $k$ for any $j = 1, \
\ldots, \ s$. The general number of the eigenvalues (taking into
account their multiplicities), equal in modulus to $\rho(A \wedge
A)$, is not greater, than $(\frac{k}{2} - 1)k$, i.e. the product
of the number of the block and the maximal index of imprimitivity
of any of them. We came to the contradiction, because there is
$C_k^2 = \frac{k(k-1)}{2}
> (\frac{k}{2} - 1)k$ eigenvalues, equal in
modulus to $\rho(A \wedge A)$.
 $\blacksquare$

 {\bf Example 1.} Let the operator $A:{\mathbb{R}}^{5} \rightarrow
{\mathbb{R}}^{5}$ be defined by the matrix $$\mathbf{A} =
\begin{pmatrix}
  0 & 1 & 0 & 0 & 0 \\
  0 & 0 & 1 & 0 & 0 \\
  0 & 0 & 0 & 1 & 0 \\
  0 & 0 & 0 & 0 & 1 \\
  1 & 0 & 0 & 0 & 0
\end{pmatrix}.$$
 This matrix is obviously nonnegative and irreducible.

In this case the second compound matrix is the following:
$${\mathbf A}^{(2)} =
\begin{pmatrix}
0 & 0 & 0 & 0 & 1 & 0 & 0 & 0 & 0 & 0 \\
0 & 0 & 0 & 0 & 0 & 1 & 0 & 0 & 0 & 0 \\
0 & 0 & 0 & 0 & 0 & 0 & 1 & 0 & 0 & 0 \\
-1 & 0 & 0 & 0 & 0 & 0 & 0 & 0 & 0 & 0 \\
0 & 0 & 0 & 0 & 0 & 0 & 0 & 1 & 0 & 0 \\
0 & 0 & 0 & 0 & 0 & 0 & 0 & 0 & 1 & 0 \\
0 & -1 & 0 & 0 & 0 & 0 & 0 & 0 & 0 & 0 \\
0 & 0 & 0 & 0 & 0 & 0 & 0 & 0 & 0 & 1 \\
0 & 0 & -1 & 0 & 0 & 0 & 0 & 0 & 0 & 0 \\
0 & 0 & 0 & -1 & 0 & 0 & 0 & 0 & 0 & 0 \\
\end{pmatrix}.$$

 It's easy to see, that the matrix ${\mathbf A}^{(2)}$ is
reducible and ${\mathcal J}$--sign-symmetric. In this case we have
$2^2 = 4$ ways of constructing the set ${\mathcal J}$:  ${\mathcal
J} = \{1, \ 5, \ 8, \ 10, \ 2, \ 6, \ 9\}$, ${\mathcal J} = \{4, \
7, \ 3\}$, ${\mathcal J} = \{1, \ 5, \ 8, \ 10, \ 7, \ 3\}$,
${\mathcal J} = \{4, \ 2, \ 6, \ 9\}$. Examine the sets $W$,
corresponding to this sets of indices ${\mathcal J}$. Such sets
$W$ defines non-transitive binary relations on the set of the
indices $\{1, \ 2, \ 3, \ 4, \ 5 \}$. The operator $A$ satisfies
the conditions of theorem 8, case (2). It's easy to see, that $A$
has the first positive simple eigenvalue $\lambda = \rho(A) = 1$,
and there is five (an odd number) eigenvalues $1$, $e^{\frac{2\pi
i}{5}}$, $e^{\frac{4\pi i}{5}}$, $e^{\frac{6\pi i}{5}}$ and
$e^{\frac{8\pi i}{5}}$ on the spectral circle $|\lambda| = 1$, all
of them are simple and coincide with $5$th roots of $1$.

Let us notice, that if we assume the irreducibility of the second
compound matrix ${\mathbf A}^{(2)}$, as well, as of the initial
matrix ${\mathbf A}$, then the theorem 8 can be refined.

 {\bf Theorem 9.} {\it Let the matrix ${\mathbf A}$ of a linear
operator $A: {\mathbb R}^n \rightarrow {\mathbb R}^n$ be
${\mathcal J}$--sign-symmetric and irreducible. Let its second
compound matrix ${\mathbf A}^{(2)}$ be also ${\mathcal
J}$--sign-symmetric and irreducible. Then one of the following two
cases takes place:
\begin{enumerate}
\item[\rm (1)] The set $\widehat{W}({\mathcal J},
\widetilde{{\mathcal J}})$ is transitive. Then $h(A) = 1$, $h(A
\wedge A)$ is an arbitrary, and the operator $A$ has two positive
simple eigenvalues $\lambda_{1}$, $\lambda_{2}$:
$$\rho(A) = \lambda_{1} > \lambda_{2} \geq |\lambda_{3}| \geq
\ldots \geq |\lambda_{n}|. $$ If $h(A) = h(A \wedge A) = 1$, then
$\lambda_{2}$ is different in modulus from the other eigenvalues.
If $h(A) = 1$, and $h(A \wedge A)
> 1$, then the operator $A$ has $h(A \wedge A)$ eigenvalues $\lambda_{2}, \lambda_{3}, \ldots, \lambda_{h(A \wedge A) + 1}
$, equal in modulus to $\lambda_{2}$, each of them is simple, and
they coincide with the $h(A \wedge A)$th roots from
$\lambda_{2}^{h(A \wedge A)}$.
\item[\rm (2)] the set $\widehat{W}({\mathcal J},
\widetilde{{\mathcal J}})$ is not transitive. Then $h(A) = h(A
\wedge A) = 3$, and there is just three eigenvalues on the
spectral circle $|\lambda| = \rho(A)$. Each of them is simple, and
they coincide with the $3$th roots of $(\rho(A))^3$.
\end{enumerate}}

The proof of the theorem 5 can be found in [3] (see [3], theorem
15).

Using lemma 1 and lemma 3, we can prove one else theorem.

{\bf Theorem 10.} {\it Let the matrix ${\mathbf A}$ of a linear
operator $A$ be ${\mathcal J}$--sign-symmetric and irreducible.
Let its second compound matrix ${\mathbf A}^{(2)}$ be ${\mathcal
J}$--sign-symmetric. Then the operator $A$ has a simple positive
eigenvalue $\lambda_{1} = \rho(A)$. Moreover, if at least one
element $a_{ii}$, situated on the principal diagonal of $\mathbf
A$, is not equal to zero, then $h(A) = 1$, and the second in
modulus eigenvalue $\lambda_{2}$ of the operator $A$ is
nonnegative and different in modulus from the first eigenvalue:
 $$0 \leq \lambda_2 < \lambda_{1}.$$ If ${\mathbf A}$ has at least
one positive principal minor of the second order, then $\lambda_2
> 0$.}

Note, that the conditions for the diagonal elements, given in
theorem 10, are only sufficient for the existence of the second in
modulus nonnegative eigenvalue of a $2$--totally ${\mathcal
J}$--sign-symmetric matrix, unlike the conditions for the set $W$,
given in theorem 9, which are necessary and sufficient.

\section{Generalization of the Gantmacher--Krein theorems to the
case of a reducible $2$--totally ${\mathcal J}$--sign-symmetric
matrix.}

Let us describe the spectrum of a reducible ${\mathcal
J}$--sign-symmetric matrix with a ${\mathcal J}$--sign-symmetric
second compound matrix, using the theorems proved above.

{\bf Theorem 11.} {\it Let the matrix ${\mathbf A}$ of a linear
operator $A$ be ${\mathcal J}$--sign-symmetric. Let its second
compound matrix ${\mathbf A}^{(2)}$  be also ${\mathcal
J}$--sign-symmetric. Let, in addition, $\rho(A) > 0$. Then the
operator $A$ has a positive eigenvalue $\lambda_{1} = \rho(A)$.
Let $m \geq 1$ be the multiplicity of the eigenvalue $\rho(A)$.
Then there is $m$ sets of eigenvalues on the largest spectral
circle $|\lambda| = \rho(A)$, with an odd number $k_j \geq 1 \ \
(j = 1, \ \ldots, \ m)$ of eigenvalues in the $j$-th set. The
eigenvalues of the $r$-th set coincide with the $k_j$-th roots of
$(\rho(A))^{k_j}$.}

$\square$ Apply theorem 6 to the matrix ${\mathbf A}$. We get,
that there are $m \geq 1$ of irreducible blocs ${\mathbf A}_j \ \
(j = 1, \ \ldots, \ m)$ with $\rho(A_j) = \rho(A)$ on the
principal diagonal of form (6) of the matrix ${\mathbf A}$. Apply
theorem 8 to every submatrix ${\mathbf A}_j$, which are, as it
follows from lemma 2, irreducible and ${\mathcal
J}$--sign-symmetric, with
 ${\mathcal J}$--sign-symmetric second compound matrices
${\mathbf A}_j^{(2)}$. We get, that there is an odd number $k_j
\geq 1$ of eigenvalues on the spectral circle $|\lambda| =
\rho(A_j)$. All of them are simple and coincide with the $k_j$th
roots of $(\rho(A))^{k_j}$. Applying theorem 6 once again, we get,
that all this eigenvalues are the eigenvalues of the operator $A$.
$\blacksquare$

\section{Concluding remarks} The results of this article can be easily generalized to the case
of $k$-totally ${\mathcal J}$-sign-symmetric matrices with $k = 3,
\ 4, \ 5, \ \ldots.$

\section*{References}

1. {\bf F.R. Gantmacher, M.G. Krein.} {\it Oscillation Matrices
and Kernels and Small Vibrations of Mechanical Systems.}
--- AMS Bookstore, 2002. --- 310 p.

2. {\bf O.Y. Kushel, P.P. Zabreiko.} {\it Gantmakher-Krein theorem
for by-non-negative operators in spaces of  functions}  //
Abstract and Applied Analysis. --- Vol. 2006. Article ID 48132.
--- P. 1-15.

3. {\bf O.Y. Kushel.} {\it On spectrum and approximations of one
class of sign-symmetric matrices.}
--- to appear.

4. {\bf Tsoy-Wo Ma.} {\it Classical analysis on normed spaces.}
--- World Scientific Publishing, 1995. --- 356 p.

5. {\bf I.M. Glazman, Yu.I. Liubich.} {\it Finite-dimensional
linear analysis.} --- M.: Nauka, 1969. --- 476 p. (Russian)

6. {\bf K. Kuratovski.} {\it Topology, I, II, revised 2nd ed.} ---
Academic Press, New York, 1966.

7. {\bf H. Minc.} {\it Nonnegative matrices.}
--- John Wiley and Sons, New York, 1988.

8. {\bf F. Gantmacher.} {\it The Theory of Matrices.} Volume 1,
Volume 2.
--- Chelsea. Publ. New York, 1990.

\medskip

O.Y. Kushel

Belorussian State University

address: 220050, Republic of Belarus, Minsk, Nezavisimosti sq., 4.

e-mail: kushel@mail.ru

\end{document}